\documentclass[11pt]{amsart}
\usepackage[utf8]{inputenc}
\usepackage{amsmath,amsfonts,amssymb,amsthm}
\usepackage{latexsym,bm,graphicx}
\usepackage{mathrsfs}
\usepackage{color}
\usepackage{hyperref}

\title[Harmonic functions with polynomial growth]{Harmonic functions with polynomial growth on manifolds with nonnegative Ricci curvature}
\author{Xian-Tao Huang}
\address{School of Mathematics\\  Sun Yat-sen University\\ Guangzhou 510275\\ E-mail address: hxiant@mail2.sysu.edu.cn}

\newtheorem{thm}{Theorem}[section]
\newtheorem{prop}[thm]{Proposition}

\newtheorem{cor}[thm]{Corollary}

\newtheorem{prob}[thm]{Problem}

\theoremstyle{definition}
\theoremstyle{remark}

\newtheorem{defn}[thm]{Definition}
\newtheorem{rem}[thm]{Remark}

\numberwithin{equation}{section}

\newcommand\tbbint{{-\mkern -16mu\int}}

\newcommand\dbbint{{-\mkern -19mu\int}}

\newcommand\bbint{
{\mathchoice{\dbbint}{\tbbint}{\tbbint}{\tbbint}}
}

\begin{document}

\maketitle
\begin{abstract} Suppose $(M,g)$ is a Riemannian manifold having dimension $n$, nonnegative Ricci curvature, maximal volume growth and unique tangent cone at infinity.
In this case, the tangent cone at infinity $C(X)$ is an Euclidean cone over the cross-section $X$.
Denote by $\alpha=\lim_{r\rightarrow\infty}\frac{\mathrm{Vol}(B_{r}(p))}{r^{n}}$ the asymptotic volume ratio.
Let $h_{k}=h_{k}(M)$ be the dimension of the space of harmonic functions with polynomial growth of growth order at most $k$.
In this paper, we prove a upper bound of $h_{k}$ in terms of the counting function of eigenvalues of $X$.
As a corollary, we obtain $\lim_{k\rightarrow\infty}k^{1-n}h_{k}=\frac{2\alpha}{(n-1)!\omega_{n}}$.
These results are sharp, as they recover the corresponding well-known properties of $h_{k}(\mathbb{R}^{n})$.
In particular, these results hold on manifolds with nonnegative sectional curvature and maximal volume growth.

\vspace*{5pt}

\vspace*{5pt}

\end{abstract}
\section{Introduction}  

Suppose $(M^{n},g)$ is a noncompact complete $n$-dimensional manifold.
Throughout this paper, we always assume $n\geq2$, and the manifold $M$ will assume to be noncompact and have nonnegative Ricci curvature if not explicitly mentioned.
We denote by $\mu$ the volume element induced by $g$.
We fix a point $p\in M$, and denote by $\rho(x)=d(x,p)$.
For a $k>0$, we consider the linear space
$$\mathcal{H}_{k}(M)=\{u\in C^{\infty}(M)\mid \Delta u=0, |u(x)|\leq C(\rho(x)^{k}+1) \text{ for some }C>0\}$$
and denote by $h_{k}(M)= \mathrm{dim}\mathcal{H}_{k}(M)$.

It is well-known that $\mathcal{H}_{k}(\mathbb{R}^{n})$ consists of harmonic polynomials of degree at most $k$, and $h_{k}(\mathbb{R}^{n})\thicksim\frac{2}{(n-1)!}k^{n-1}$ as $k\rightarrow\infty$.
See Appendix B of \cite{L97} for the details of these facts.

There are many researches on polynomial growth harmonic functions on manifolds with nonnegative Ricci curvature.

Yau \cite{Y75} proved that, on complete manifolds with nonnegative Ricci curvature, any positive harmonic function is constant.
In \cite{Cg82}, Cheng further proved that on such manifolds any harmonic function of sublinear growth must be constant.

A famous conjecture of Yau says that on a manifold $(M^{n},g)$ with non-negative Ricci curvature it always holds $h_{k}(M)<\infty$.
Li and Tam solved the case for $k=1$ and the case $n=2$ of this conjecture (see \cite{LT89} \cite{LT91}).
Yau's conjecture was completely solved by Colding and Minicozzi in \cite{CM97a}.
Later on, there are many researches which give better estimates of $h_{k}(M)$.

For example, a more precise upper bound for the dimensions was obtained:

\begin{thm}[see \cite{CM98b}, \cite{L97}]\label{upperbound}
If $(M^{n},g)$ has nonnegative Ricci curvature, then $h_{k}(M)\leq Ck^{n-1}$ for all $k\geq1$, where $C=C(n)>0$ is a constant depending only on $n$.
\end{thm}
Note that the power $n-1$ in Theorem \ref{upperbound} is sharp compared to the Euclidean case.

There was once a stronger conjecture saying that if $(M^{n},g)$ has nonnegative Ricci curvature, then $h_{k}(M^{n})\leq h_{k}(\mathbb{R}^{n})$ holds for every $k$.
In 2001, Donnelly (see \cite{Don01}) gave a counterexample to this conjecture.
However, there is still an open problem (see e.g. the last question in Section 28 of \cite{L12}) on the upper bound of dimensions as follows.
\begin{prob}\label{ques1.1}
If $(M^{n},g)$ has nonnegative sectional curvature, then is it true that
\begin{align}\label{1.0001}
h_{k}(M^{n})\leq h_{k}(\mathbb{R}^{n})
\end{align}
holds for every $k$?
\end{prob}

This problem is unsolved, except the case of $k=1$ (see \cite{LT89} and \cite{CCM}).

One may also ask whether there is a better estimate on the constant $C$ in the statement of Theorem \ref{upperbound}, at least under some geometric assumptions on $M$.

In \cite{CM98b}, Colding and Minicozzi proved that if $(M^{n},g)$ has nonnegative Ricci curvature then
\begin{align}\label{cm_asy}
h_{k}(M)\leq C_{1} \alpha k^{n-1}+C_{2}f(k^{n-1}),
\end{align}
where
\begin{align}\label{1.01111}
\alpha:=\lim_{r\rightarrow\infty}\frac{\mu(B_{r}(p))}{r^{n}}
\end{align}
is the  asymptotic volume ratio; $C_{1}$ and $C_{2}$ are positive constants depending only on $n$; the function $f:\mathbb{R}^{+}\rightarrow\mathbb{R}^{+}$ also depends only on $n$ and satisfies
$f(t)\leq t$ and $\lim_{t\rightarrow\infty}\frac{f(t)}{t}=0$.
See Theorem 0.26 in \cite{CM98b}.

In \cite{Huang20}, the author obtained the following theorem.

\begin{thm}[see \cite{Huang20}]\label{main-4}
Suppose $(M^{n}, g)$ is a complete manifold with nonnegative Ricci curvature and maximal volume growth, i.e. $\alpha>0$.
Assume the tangent cone at infinity of $M$ is unique.
Then we have
\begin{align}\label{1.22221}
\lim_{k\rightarrow\infty} k^{-n}\sum_{i=1}^{k}h_{i-1}=\frac{2\alpha}{n!\omega_{n}},
\end{align}
and
\begin{align}\label{1.22222}
\liminf_{k\rightarrow\infty}k^{1-n}h_{k}=\frac{2\alpha}{(n-1)!\omega_{n}},
\end{align}
where $\omega_{n}$ is the volume of a unit ball in $\mathbb{R}^{n}$.
\end{thm}

\begin{rem}
Suppose $(M^{n}, g)$ has nonnegative Ricci curvature and is collasped, i.e. $\alpha=0$, then by (\ref{cm_asy}), we have $\lim_{k\rightarrow\infty}k^{1-n}h_{k}=0$.
In particular, (\ref{1.22221}) and (\ref{1.22222}) still hold.
\end{rem}

In 1999, Li and Wang have proved a similar theorem, where they assume the manifold having nonnegative sectional curvature, see Theorem 2.2 in \cite{LW99}.
Recall that if a manifold has nonnegative sectional curvature, then its tangent cone at infinity is a unique metric cone.
And there are many examples of manifolds satisfying the assumptions in Theorem \ref{main-4} but do not have nonnegative sectional curvature.
Also note that we obtain equalities in (\ref{1.22221}) and (\ref{1.22222}), while in \cite{LW99}, the conclusions are some inequalities.

The proof of Theorem \ref{main-4} combines Li and Wang's proof (\cite{LW99}), Cheeger-Colding's theory (see e.g. \cite{CC96}, \cite{CC97}), some recent development in \textmd{RCD} theory, as well as a theorem on the lower bound of $h_{k}(M)$ (obtained in \cite{D04}, \cite{Huang19}, see Theorem \ref{main-2} below).

A key viewpoint in the proof of Theorem \ref{main-4} is that, for a manifold $(M^{n},g)$ with nonnegative Ricci curvature and maximal volume growth, its asymptotically conic property makes the harmonic functions with polynomial growth on it behave like harmonic functions on a cone.
This viewpoint has influenced many previous works, see \cite{CM97b}, \cite{CM98b}, \cite{D04}, \cite{X15}, \cite{H15}, \cite{Huang19} etc.

A natural question is that, under the assumption of Theorem \ref{main-4}, can we improve the $\liminf$ in (\ref{1.22222}) to $\lim$ ?
The main results of this paper give an affirmative answer to this question.
In fact, we will prove some stronger results, see Theorem \ref{main-6} below.

Before we state our main results, we review some work concerning the lower bound of $h_{k}(M)$.

There are manifolds with nonnegative Ricci curvature on which any harmonic function with polynomial growth is constant, see Theorem 1 in \cite{S00} and Corollary 8.12 of \cite{CM98b}.
In other words, to ensure the existence of nontrivial harmonic functions with polynomial growth, some additional assumptions on the manifold are necessary.

In Theorem 0.1 of \cite{D04}, Ding proved that suppose $(M^{n},g)$ has nonnegative Ricci curvature, maximal volume growth and has a unique tangent cone at infinity, then $h_{k}(M)\geq C k^{n-1}$ for some constant $C$ depending on $n$ and $\alpha$.
Note that the power $n-1$ in this result is sharp compared to the Euclidean case.
In a recent paper, Xu obtained the existence of nontrivial harmonic functions with polynomial growth on some manifolds satisfying certain assumptions on its tangent cone at infinity, see Theorem 1.7 in \cite{X15}.
Note that Xu's existence theorem does not require the manifold to have maximal volume growth, and an explicit example of a collapsed manifold satisfying the assumptions in Xu's existence theorem was constructed in \cite{X15}.

In \cite{D04} and \cite{X15}, the ideas to find nontrivial harmonic functions with polynomial growth are similar.
Making use of the conic structure at infinity, we transplant suitable harmonic functions on the tangent cone at infinity back to obtain harmonic functions defined on larger and larger geodesic balls of the manifold.
By construction and making use of certain three circles theorems, we can prove these harmonic functions satisfy the desired growth rate.
Finally by Arzela-Ascoli theorem we construct a harmonic function with the desired growth rate on the whole  manifold.
In \cite{X15}, Xu proved a three circles theorem, which is different from the one given by Ding, see Theorem 3.2 in \cite{X15} and Lemma 1.1 in \cite{D04} respectively.

Later on, following the ideas of Xu, the author obtained an improved lower bound of $h_{k}(M)$ when $M$ satisfies the same assumptions as in Xu's existence theorem, see Theorem 1.2 in \cite{Huang19}.
This result recovers Theorem 0.1 in \cite{D04} if the manifold has maximal volume growth.
In the following, we state the lower bound theorem in the special case that the manifold has maximal volume growth, the more general cases can consult Theorem 1.2 in \cite{Huang19}.

\begin{thm}[See \cite{D04} and \cite{Huang19}]\label{main-2}
Suppose $(M^{n},g)$ is a complete Riemannian manifold with nonnegative Ricci curvature with maximal volume growth.
Assume $M$ has a unique tangent cone at infinity, denoted by $(C(X),d_{C(X)},m_{C(X)})$, which is an Euclidean cone over the cross section $(X, d_{X}, m_{X})$.
Let $0=\lambda_{0}<\lambda_{1}\leq\lambda_{2}\leq\ldots$ be the eigenvalues (counted with multiplicity) for the Laplacian operator on $(X, d_{X}, m_{X})$.
Let $N_{X}:\mathbb{R}^{+}\rightarrow \mathbb{Z}^{+}$ be the counting function with respect to $(X,d_{X},m_{X})$, i.e.
\begin{align}\label{countfcn}
N_{X}(\lambda) := \#\{i \in \mathbb{Z}^{+}\cup\{0\} |\lambda_{i}\leq\lambda\}.
\end{align}
Then given any $k>0$ such that $k(k+n-2)>\lambda_{1}$,
\begin{align}\label{5.33333}
h_{k}\geq N_{X}(\beta(\beta+n-2))
\end{align}
holds for any positive number $\beta<k$.
\end{thm}

The key point in the proof of Theorem \ref{main-2} is based on the method in \cite{X15}: we use $L^{2}$-orthonormal eigenfunctions of the cross section $X$ to construct as many as possible linear independent harmonic functions  with polynormial growth on $M$.

Finally we note that in a more recent paper \cite{X20}, Xu proved another existence theorem provided the manifold satisfies a quantitative strong unique continuation in the sense of Colding and Minicozzi (see \cite{CM97c}).

Now we begin to introduce the new results in this paper.

\begin{thm}\label{main-5}
Suppose $(M^{n},g)$ is a complete Riemannian manifold with nonnegative Ricci curvature with maximal volume growth.
Assume $M$ has a unique tangent cone at infinity, denoted by $(C(X),d_{C(X)},m_{C(X)})$, which is an Euclidean cone over the cross section $(X, d_{X}, m_{X})$.
Let $0=\lambda_{0}<\lambda_{1}\leq\lambda_{2}\leq\ldots$ be the eigenvalues (counted with multiplicity) for the Laplacian operator on $(X, d_{X}, m_{X})$.
Let $N_{X}:\mathbb{R}^{+}\rightarrow \mathbb{Z}^{+}$ be the counting function with respect to $(X,d_{X},m_{X})$ as in (\ref{countfcn}).
Then given any $k>0$, we have
\begin{align}\label{5.33333123}
h_{k}\leq N_{X}(k(k+n-2)).
\end{align}
\end{thm}

\begin{rem}
If we take $k>0$ satisfying $k(k+n-2)<\lambda_{1}$, then by (\ref{5.33333123}), we have $h_{k}=1$.
This is a special case of Liouville type theorems proved by Honda (see \cite{H15}) and the author (see \cite{Huang19}), where the conclusions hold without the assumption on the uniqueness of tangent cone at infinity.
Note that the proof of Theorem 1.4 in \cite{Huang19} is different from the one given in \cite{H15}, which depending on the maximal volume growth condition.
Hence Theorem 1.4 in \cite{Huang19} applies to certain manifolds which does not have maximal volume growth.
Motivated by \cite{CM97b}, we define some functionals related to a harmonic function on $C(X)$ in \cite{Huang19}, and obtain some good properties of these functionals.
In the proof of Theorem 1.4 in \cite{Huang19}, given $u\in \mathcal{H}_{k}(M)$, we choose parameters carefully to blow down the manifold so that we obtain a harmonic function on $C(X)$.
By the choice of parameters and the good properties of the functionals defined above, we obtain many useful information on $u$.
\end{rem}

A key in the proof of Theorem \ref{main-5} is a three circles type theorem, see Proposition \ref{main-8-weak}.
The proof of Proposition \ref{main-8-weak} is based on a modification of that of Theorem 1.4 in \cite{Huang19}, where we will show any $u\in\mathcal{H}_{k}(M)$ has the desired three circles type property near infinity, and then we apply Xu's three circles theorem (i.e. Theorem 3.2 in \cite{X15}) to $u$.

By Theorems \ref{main-2} and \ref{main-5}, we have the following corollary.

\begin{cor}\label{main-6}
Suppose $(M^{n},g)$ satisfies the assumptions in Theorem \ref{main-5}, with the unique tangent cone at infinity, denoted by $(C(X),d_{C(X)},m_{C(X)})$.
Let $N_{X}:\mathbb{R}^{+}\rightarrow \mathbb{Z}^{+}$ be the counting function with respect to $(X,d_{X},m_{X})$ as in (\ref{countfcn}).
Denote by
\begin{align}\label{Dx}
\mathscr{D}_{X}:=\{\beta\geq 0\mid  \beta(\beta+n-2)=\lambda_{i}\text{ for some eigenvalue $\lambda_{i}$ of $(X, d_{X}, m_{X})$} \}.
\end{align}
Then given any $k>0$ such that $k\notin \mathscr{D}_{X}$, we have
\begin{align}\label{5.33333123ggggg}
h_{k}= N_{X}(k(k+n-2)).
\end{align}
\end{cor}

Combining Corollary \ref{main-6} with Weyl's law on $\mathrm{RCD}$ spaces (see \cite{AHT17} \cite{ZZ17}), we have the following corollary.

\begin{cor}\label{main-7}
Suppose $(M^{n},g)$ is a complete Riemannian manifold with nonnegative Ricci curvature and asymptotic volume ratio $\alpha>0$.
Assume that $M$ has a unique tangent cone at infinity.
Then
\begin{align}\label{1.22222345}
\lim_{k\rightarrow\infty}k^{1-n}h_{k}=\frac{2\alpha}{(n-1)!\omega_{n}}.
\end{align}
\end{cor}

Making use of Faulhaber's formula, one can check that (\ref{1.22222345}) implies (\ref{1.22221}).
Thus Corollary \ref{main-7} improves Theorem \ref{main-4}.

We also prove the following theorem.
The conclusion (\ref{5.333331234}) is a three circles theorem, which improves Proposition \ref{main-8-weak}.

\begin{thm}\label{main-8}
Suppose $(M^{n},g)$ satisfies the assumptions in Theorem \ref{main-5}, with the unique tangent cone at infinity, denoted by $(C(X),d_{C(X)},m_{C(X)})$.
Let $\mathscr{D}_{X}$ be given by (\ref{Dx}).
Suppose $k>0$ and $k(k+n-2)\geq\lambda_{1}$.
Then for any $u \in\mathcal{H}_{k}(M)$ with $u(p)=0$, there exists a $\gamma\in \mathscr{D}_{X}$ satisfying $\gamma\leq k$ and the followings.
\begin{description}
  \item[(1)] Given any $\epsilon>0$, there exists a $T=T(M,u,\epsilon)$ such that
  \begin{align}\label{5.333331234}
  \bbint_{B_{2r}(p)}u^{2}d\mu \leq 2^{2\gamma+2\epsilon}\bbint_{B_{r}(p)}u^{2} d\mu
  \end{align}
  holds for any $r\geq T$.
  \item[(2)] $u\in\mathcal{H}_{\gamma+\epsilon}(M)$ for any $\epsilon>0$.
  \item[(3)] $u\notin\mathcal{H}_{\gamma-\epsilon}(M)$ for any $\epsilon>0$.
\end{description}

\end{thm}

\begin{rem}
It is unclear whether (2) in Theorem \ref{main-8} can be improved to $u\in\mathcal{H}_{\gamma}(M)$.
\end{rem}

\begin{rem}
We give some remarks on the connection between Theorem \ref{main-5} and the open problem \ref{ques1.1}.
Suppose $(M^{n},g)$ has nonnegative sectional curvature and maximal volume growth, then its tangent cone at infinity is unique.
Since the tangent cone at infinity $(C(X), d_{C(X)})$ is an $n$-dimensional Alexandrov space with nonnegative curvature, it is well known that the cross-section $(X,d_{X})$ is an $(n-1)$-dimensional Alexandrov space with curvature bounded from below by $1$.
By Theorem \ref{main-5}, it is not hard to see that, when the manifolds has maximal volume growth, (\ref{1.0001}) is true provided we can give an affirmative answer to the following long standing open problem:
\begin{prob}
Suppose $(X,d_{X})$ is an $(n-1)$-dimensional Alexandrov space with curvature bounded from below by $1$, then is it true that
\begin{align}
\lambda_{i}(X)\geq \lambda_{i}(S^{n-1})
\end{align}
holds for every $i\geq0$?
\end{prob}
However, as far as the author knows, except the cases $i=0,1$, there is almost no conclusions on this conjecture, even in the case that $X$ is a smooth manifold.

\end{rem}

\begin{rem}\label{rem1.14}
The conclusions similar to Theorem \ref{main-5}, Corollary \ref{main-6}, and Theorem \ref{main-8} still hold provided the manifold $(M^{n}, g)$ has nonnegative Ricci curvature and satisfies the followings:
\begin{description}
  \item[(1)] its tangent cone at infinity with renormalized limit measure is a unique metric cone $C(X)$ with the unique conic measure of power $\kappa\geq2$ (see \cite{X15} for definitions of these terminologies);
  \item[(2)] $\mathcal{H}^{1}(X)>0$;
  \item[(3)] there is a positive constant $A$ such that
  \begin{align} \label{volume-condition}
  \limsup_{r\rightarrow\infty}\frac{\mu(B_{r}(p))}{r^{\kappa}}\leq A<\infty.
  \end{align}
\end{description}
The proofs are almost the same as the ones in this paper.
Note that Xu's example (see Example 4.9 in \cite{X15}) is a noncompact manifold satisfying assumptions (1)-(3) and having positive sectional curvature but does not have maximal volume growth.

The conclusions under the assumptions (1)-(3) are related to Question 6.53 and Theorem 6.58 in \cite{CM98b}.
We will discuss the relation between them in Section \ref{sec-5}.
In Section \ref{sec-5} we will also prove an asymptotic dimension estimate similar to (\ref{1.22222345}) on collapsed manifolds with nonnegative sectional curvature and satisfying an assumption weaker than assumptions (1)-(3).
\end{rem}

\begin{rem}
The methods in this paper can be applied to prove similar dimension estimates of harmonic functions with polynomial growth on certain $\mathrm{RCD}(0,N)$ spaces.
For example, conclusions similar to Theorems \ref{main-5}, \ref{main-8} and Corollaries \ref{main-6}, \ref{main-7} hold on any non-collapsed (in the sense of \cite{DePGil18}) $\mathrm{RCD}(0,N)$ spaces $(Y,d,\mathcal{H}^{N})$ with $\alpha:=\lim_{r\rightarrow\infty}\frac{\mathcal{H}^{N}(B_{r}(p))}{r^{N}}>0$ for some $p\in Y$ and whose tangent cone at infinity is unique.
Here $N$ is a positive integer and $\mathcal{H}^{N}$ is the $N$-dimensional Hausdorff measure.
\end{rem}


\vspace*{20pt}

\noindent\textbf{Acknowledgments.}
The author would like to thank Prof. B-L Chen, H-C Zhang, X-P Zhu and Dr. H-Z Huang for discussions.
The author is partially supported by NSFC (Nos. 12025109 and 11521101).

\section{Preliminaries and notations}\label{sec-2}

In this section, we recall some background concepts on calculus and properties of Ricci limit spaces.

Let $(Z,d,\nu)$ be a complete path-connected metric measure space, equipped with a Radon measure $\nu$.
For a locally Lipschitz function $f:Z\rightarrow \mathbb{R}$, denote the pointwise Lipschitz constant of $f$ by
$$\textmd{Lip} f(x)=\limsup_{d(z,x)\rightarrow0}\frac{|f(z)-f(x)|}{d(z,x)}.$$

For $f\in L^{2}(Z)$, the Cheeger energy is defined by
$$\textmd{Ch}(f)=\inf_{\{f_{n}\}}\biggl\{\liminf_{n\rightarrow\infty}\frac{1}{2}\int(\textmd{Lip}(f_{n}))^{2}d\nu \bigl| f_{n}\in \textmd{Lip}(Z)\cap L^{2}(Z),\| f_{n}-f \|_{L^{2}}\rightarrow0\biggr\}.$$

The Sobolev space $W^{1,2}(Z)$ is defined as $W^{1,2}(Z):= \bigl\{f\in L^{2}(Z)\bigl|\mathrm{Ch}(f)< \infty\bigr\}$.
$W^{1,2}(Z)$ is equipped with the norm
$$\| f\|^{2}_{W^{1,2}}:=\| f\|^{2}_{L^{2}}+2\mathrm{Ch}(f).$$
It is known that for any $f\in W^{1,2}(Z)$, there exists $|Df|\in L^{2}(Z)$ such that $2\mathrm{Ch}(f)=\int_{Z}|Df|^{2}d\nu$.
$|Df|$ is called the minimal weak upper gradient of $f$.
The minimal weak upper gradient is local in the sense that for any $f,g\in W^{1,2}(Z)$, $|Df|=|Dg|$ holds $\nu$-a.e. on $\{x\in Z\bigl| f(x)=g(x)\}$.

For an open set $U\subset Z$, we define $W^{1,2}(U)$ to be the space of functions $f:U\rightarrow\mathbb{R}$ locally equal to some function in $W^{1,2}(Z)$ and satisfy $f, |Df|\in L^{2}(U)$.
We use $W^{1,2}_{\mathrm{loc}}(U)$ to denote the space of functions $f:U\rightarrow\mathbb{R}$ locally equal to some function in $W^{1,2}(Z)$.

Let $U\subset X$ be an open set, then for any $f,g\in W^{1,2}_{\mathrm{loc}}(U)$, $\langle D f,D g\rangle:U\rightarrow\mathbb{R}$ is $\nu$-a.e. defined to be
$$\langle D f,D g\rangle:=\inf_{\epsilon>0}\frac{|D(g+\epsilon f)|^{2}-|Dg|^{2}}{2\epsilon},$$
where the infimum is in $\nu$-essential sense.
In the case that $W^{1,2}(Z)$ is a Hilbert space (this holds if $(Z, d, \nu)$ is a Ricci-limit space), the map $W^{1,2}_{\mathrm{loc}}(U)\ni f, g\mapsto \langle D f,D g\rangle\in L^{1}_{\mathrm{loc}}(U)$ is bilinear and symmetric, and we have $\langle D f,D f\rangle=|Df|^{2}$.

Given an open set $U\subset Z$, $D({\Delta},U)\subset W^{1,2}_{\mathrm{loc}}(U)$ is the space of $f\in W^{1,2}_{\mathrm{loc}}(U)$ such that there exists a signed Radon measures $\mu$ on $U$ such that
$$\int g d\mu =-\int \langle D f,D g\rangle d\nu$$
holds for any $g: Z\rightarrow\mathbb{R}$ Lipschitz with $\mathrm{supp}(g)\subset\subset U$.
$\mu$ is uniquely determined and we denote it by ${\Delta}f$.
If $U=Z$, $f\in W^{1,2}(Z)\cap D({\Delta},Z)$ and ${\Delta}f=hm$ for some $h\in L^{2}(Z)$, then we say $f\in D(\Delta)$, and denote by $\Delta f=h$.

A function $f\in W^{1,2}_{\mathrm{loc}}(U)$ is called harmonic on $U$ if $u\in D({\Delta},U)$ and ${\Delta} u=0$ in $U$.

For more information on calculus in general metric measure spaces, the readers can refer to \cite{AGS14-1} \cite{C99} \cite{Gig15} etc.

In the following, we assume $(M^{n},g)$ is a complete Riemannian manifold with nonnegative Ricci curvature.
Let $\rho$ be the distance determined by $g$, $\mu$ be the volume measure determined by $g$.
Given $r_{i}\rightarrow\infty$, we consider the rescaled metric $g_{i}=r_{i}^{-2}g$.
We denote $(M_{i}, p_{i}, \rho_{i},\nu_{i})$, where $M_{i}$ is the same differential manifold as $M^{n}$, $p_{i}=p$, $\rho_{i}$ is the distance determined by $g_{i}$, and $\nu_{i}$ is the renormalized measure defined by
$$\nu_{i}(A):=\frac{1}{\mu_{i}(B_{1}(p_{i}))}\mu_{i}(A),$$
where $A\subset M_{i}$, and $\mu_{i}$ is the volume measure determined by $g_{i}$.
Then by Gromov's compactness theorem and Theorem 1.6 in \cite{CC97}, after extracting a subsequence, we can assume $(M_{i}, p_{i}, \rho_{i},\nu_{i})$ converges to some $(M_{\infty}, p_{\infty}, \rho_{\infty},\nu_{\infty})$ in the pointed measure Gromov-Hausdorff sense.
$(M_{\infty},p_{\infty},\rho_{\infty})$ is called a tangent cone at infinity of $M$.
In general, the tangent cones at infinity may not be unique, that is, it depends on the choice of $\{r_{i}\}$ and the subsequence.
Similarly, the limit measure $\nu_{\infty}$ may depend on the choice of $\{r_{i}\}$ and the subsequence.

If $(M^{n},g)$ has maximal volume growth, i.e.
\begin{align}\label{1.223}
\alpha:=\lim_{r\rightarrow\infty}\frac{\mu(B_{r}(p))}{r^{n}}>0,
\end{align}
then by the Cheeger-Colding theory (see \cite{CC96} \cite{CC97} \cite{CC00II} etc.) and $\mathrm{RCD}$ theory (see \cite{AGMR00} \cite{EKS15} \cite{Gig13} etc.), $(M_{\infty}, p_{\infty}, \rho_{\infty},\nu_{\infty})$ has lots of nice properties.

Firstly, $(M_{\infty}, \rho_{\infty})$ is a Euclidean cone over a metric measure $(X,d_{X})$ with $\textmd{diam}(X)\leq\pi$, which is denoted by $(C(X),d_{C(X)})$, with $p_{\infty}$ being the tips.
In addition, in this case the limit measure is unique.
We have $\nu_{\infty}=\frac{1}{\alpha}\mathcal{H}^{n}$, where $\mathcal{H}^{n}$ is the $n$-dimensional Hausdorff measure.

In order to emphasis the cone structure, in the remaining part of this paper, we use the notation $(C(X),p_{\infty},d_{C(X)},m_{C(X)})$ instead of $(M_{\infty},p_{\infty},\rho_{\infty}, \nu_{\infty})$.

Secondly, the cross section $X$ has an induced metric measure structure $(X,d_{X},m_{X})$ with $m_{X}=\frac{1}{\alpha}\mathcal{H}^{n-1}$, and for any measurable $\Omega\subset\subset C(X)$, we have
\begin{align}\label{conic3}
m_{C(X)}(\Omega)=\int_{0}^{\infty}s^{n-1}ds\int_{X}\chi(\Omega_{s})dm_{X},
\end{align}
where $\Omega_{s}=\{x\in X\mid z=(x,s)\in\Omega\}$, and $\chi(\cdot)$ is the characteristic function on $X$.
It is not hard to check $\mathcal{H}^{n-1}(X)=n\alpha$.

Thirdly, $(C(X), d_{C(X)},m_{C(X)})$ is an $\textmd{RCD}(0, n)$ space, and by Corollary 1.3 of \cite{Ke15}, $(X, d_{X}, m_{X})$ is an $\textmd{RCD}(n-2,n-1)$ space.

Let $\Delta_{X}$ be the Laplacian operator on $(X,d_{X},m_{X})$.
Let $0=\lambda_{0}<\lambda_{1}\leq\lambda_{2}\leq\ldots$  be the eigenvalues (counted with multiplicity) of $\Delta_{X}$, $\{\varphi_{i}\}_{i=0}^{\infty}$ be the corresponding eigenfunctions, i.e. $\varphi_{i}\in W^{1,2}(X)\cap D({\Delta}_{X})$ and
\begin{align}\label{2.0000}
-\Delta_{X}\varphi_{i}(x)=\lambda_{i}\varphi_{i}(x).
\end{align}
Every $\varphi_{i}$ always has a Lipschitz representative in the corresponding Sobolev class (see \cite{J14}).
In the remaining part of this paper, the $\varphi_{i}$'s are always required to be Lipschitz.
We also require
\begin{align}\label{2.0001}
\int_{X}|\varphi_{i}|^{2}dm_{X}=1\text{\qquad and \qquad}\int_{X}\varphi_{i}\varphi_{j}dm_{X}=0
\end{align}
for every $i\neq j$.
By a standard Rellich type compactness Theorem argument, we can derive that $\lambda_{i}\rightarrow\infty$ and $\{\varphi_{i}\}_{i=0}^{\infty}$ spans $L^{2}(X)$.
Let $N_{X}:\mathbb{R}^{+}\rightarrow \mathbb{Z}^{+}$ be the counting function defined in (\ref{countfcn}).

Recently, Weyl's law for eigenvalues of the Laplacian operator on compact $\textmd{RCD}(K,N)$ spaces has been proved (see \cite{AHT17} \cite{ZZ17}).
Applying Corollary 4.4 in \cite{AHT17} to the $\textmd{RCD}(n-2,n-1)$ space $(X, d_{X}, m_{X})$, we can prove the following result.
See Proposition 2.6 in \cite{Huang20} for details.

\begin{prop}[see \cite{Huang20}]\label{prop2.1111}
Suppose $(C(X),d_{C(X)},m_{C(X)})$ is a tangent cone at infinity of a manifold $(M^{n},g)$ with nonnegative Ricci curvature satisfying (\ref{1.223}).
Then the following Weyl's law holds:
\begin{align}\label{1.11117}
\lim_{\lambda\rightarrow\infty}\frac{N_{X}(\lambda)}{\lambda^{\frac{n-1}{2}}}=\frac{n\omega_{n-1}\alpha}{(2\pi)^{n-1}}.
\end{align}
Here $\omega_{n-1}$ is the volume of a unit ball in the $(n-1)$-dimensional Euclidean space
\end{prop}

We recall some notions and properties on convergence of functions defined on varying metric measure spaces.

\begin{defn}
Suppose metric spaces $(X_{i}, d_{i})$ pointed Gromov-Hausdorff converge to $(X_{\infty}, d_{\infty})$, with a sequence of $\epsilon_{i}$-Gromov-Hausdorff approximations $\Phi_{i}: X_{i}\rightarrow X_{\infty}$, where $\epsilon_{i}\rightarrow0$.
Suppose $f_{i}$ is a function on $X_{i}$ and $f_{\infty}$ is a function on $X_{\infty}$.
Suppose $K\subset X_{\infty}$.
If for every $\epsilon>0$, there exists $\delta>0$ such that $|f_{i}(x_{i})-f_{\infty}(x_{\infty})|<\epsilon$ holds for every $i\geq\delta^{-1}$, $x_{i}\in X_{i}$, $x_{\infty}\in K$ with $d_{\infty}(\Phi_{i}(x_{i}),x_{\infty})< \delta$, then we say $f_{i}$ converge to $f_{\infty}$ uniformly on $K$.
\end{defn}

The following theorem is a generalization of the classical Arzela-Ascoli Theorem, see e.g. Proposition 27.20 in \cite{Vi09}.
\begin{prop}\label{AA}
Suppose $(X_{i}, p_{i}, d_{i})$ pointed Gromov-Hausdorff converge to $(X_{\infty}, p_{\infty}, d_{\infty})$.
Let $R\in(0,\infty]$.
Suppose for every $i$, $f_{i}$ is a Lipschitz function defined on $B_{R}(p_{i})\subset X_{i}$ and $\mathrm{Lip}f_{i}\leq L$ on $B_{R}(p_{i})$, $| f_{i}(p_{i})|\leq C$  for some uniform constants $L$ and $C$.
Then there exits a subsequence of $f_{i}$, still denoted by $f_{i}$, and a Lipschitz function $f_{\infty}: B_{R}(p_{\infty})\rightarrow\mathbb{R}$ with $\mathrm{Lip}f_{\infty}\leq L$ such that $f_{i}$ converge uniformly to $f_{\infty}$.
\end{prop}

We also need the following useful result on convergence of harmonic functions defined on varying spaces.

\begin{prop}[\cite{H11}, see also \cite{D02} and \cite{Xu14}]\label{2.7777}
Let $(M_{i}, p_{i}, \rho_{i},\nu_{i})$ be metric measure structure induced by $(M_{i}^{n},g_{i})$ with $\mathrm{Ric}_{g_{i}}\geq K$.
Suppose $(M_{i}, p_{i}, \rho_{i},\nu_{i})$ converge in the pointed measure Gromov-Hausdorff sense to $(M_{\infty}, p_{\infty}, \rho_{\infty}, \nu_{\infty})$.
Let $f_{i}$, $g_{i}$ be harmonic functions defined on $B_{R}(p_{i})\subset M_{i}$ such that $\mathrm{Lip}f_{i}\leq L$, $\mathrm{Lip}g_{i}\leq L$ for every $i$ and some uniform constant $L$.
Suppose $f_{i}\rightarrow f_{\infty}$ and $g_{i}\rightarrow g_{\infty}$ uniformly on $B_{R}(p_{\infty})\subset M_{\infty}$, then $f_{\infty}$, $g_{\infty}$ are harmonic on $B_{R}(p_{\infty})$, and
\begin{align}
\lim_{i\rightarrow\infty}\int_{B_{r}(p_{i})}f_{i}^{2}d\nu_{i}
=\int_{B_{r}(p_{\infty})}f_{\infty}^{2}d\nu_{\infty},
\end{align}
\begin{align}
\lim_{i\rightarrow\infty}\int_{B_{r}(p_{i})}f_{i}g_{i}d\nu_{i}
=\int_{B_{r}(p_{\infty})}f_{\infty}g_{\infty}d\nu_{\infty},
\end{align}
\begin{align}
\lim_{i\rightarrow\infty}\int_{B_{r}(p_{i})}|\nabla^{(i)}f_{i}|^{2}d\nu_{i}
=\int_{B_{r}(p_{\infty})}|Df_{\infty}|^{2}d\nu_{\infty},
\end{align}
\begin{align}
\lim_{i\rightarrow\infty}\int_{B_{r}(p_{i})}\langle\nabla^{(i)}f_{i},\nabla^{(i)}g_{i}\rangle d\nu_{i}
=\int_{B_{r}(p_{\infty})}\langle Df_{\infty},Dg_{\infty}\rangle d\nu_{\infty}
\end{align}
hold for any $r\in(0,R)$, where $\nabla^{(i)}$ is the gradient with respect to $g_{i}$, and $D$ is the Cheeger derivative with respect to $(M_{\infty}, \rho_{\infty}, \nu_{\infty})$.
\end{prop}

We remark that there is a general theory on convergence of functions defined on varying $\mathrm{RCD}(K,N)$ spaces.
Except for the notion of uniform convergence, there are other notions of convergence such as (locally) $L^{2}$-weak convergence, (locally) $L^{2}$-strong convergence, (locally) $W^{1,2}$-convergence, see e.g. \cite{AH17}.
Theorem \ref{2.7777} can be generalized to the $\mathrm{RCD}(K,N)$ setting, see e.g. Theorem 4.4 and Corollary 4.5 of \cite{AH17}.

\section{Harmonic functions on conic Ricci limit space}\label{sec-3}

In the following of this paper, we assume $(M^{n},g)$ is a complete manifold having nonnegative Ricci curvature and maximal volume growth.
Let $(C(X),p_{\infty},d_{C(X)},m_{C(X)})$ be one of its tangent cones at infinity.
In the remaining part of this paper, when we talk about a harmonic function defined on $U\subset C(X)$, we always mean the Lipschitz representative (\cite{J14}) in the corresponding Sobolev class.

Let $0=\lambda_{0}<\lambda_{1}\leq\lambda_{2}\leq\ldots$ be the eigenvalues (counted with multiplicity) of the Laplacian $\Delta_{X}$, $\{\varphi_{i}(x)\}_{i=1}^{\infty}$ be the functions satisfying (\ref{2.0000}) (\ref{2.0001}).
Note that $\varphi_{i}$ is always required to be Lipschitz.
Let $N_{X}:\mathbb{R}^{+}\rightarrow \mathbb{Z}^{+}$ be the counting function defined in (\ref{countfcn}).
Suppose $\mathscr{D}_{X}$ is given by (\ref{Dx}).

It is not hard to check that if $u$ is a function on $C(X)$ given by $u(x,r)=r^{\alpha_{i}}\varphi_{i}(x)$,
where $-\Delta_{X}\varphi_{i}(x)=\lambda_{i}\varphi_{i}(x)$, $\alpha_{i}\geq0$ and $\lambda_{i}=\alpha_{i}(\alpha_{i}+n-2)$,
then $u$ is harmonic.

The following theorem, which is well-known to experts, have appeared in many papers, see \cite{H15}, \cite{X15}, \cite{Huang19}, \cite{Huang20}, \cite{CJN21} etc.
The readers can refer to \cite{Huang20} or \cite{CJN21} for a detailed proof.

\begin{thm}\label{eigenfcn}
If $u$ is a harmonic function defined on $B_{R}(p_{\infty})\subset C(X)$ with $u(p_{\infty})=0$,
then
\begin{align}\label{har_fcn}
u(x,r)=\sum_{i=1}^{\infty}c_{i}r^{\alpha_{i}}\varphi_{i}(x),
\end{align}
where the convergence in (\ref{har_fcn}) is locally uniformly on $B_{R}(p_{\infty})$ and in $W^{1,2}_{\mathrm{loc}}(B_{R}(p_{\infty}))$ sense.
In (\ref{har_fcn}), $c_{i}$ are constants, $\lambda_{i}$ and $\varphi_{i}$ satisfy (\ref{2.0000}) (\ref{2.0001}), and $\lambda_{i}=\alpha_{i}(\alpha_{i}+n-2)$, $\alpha_{i}>0$.
\end{thm}

Given a non-constant harmonic function $u(x,r)$ defined on $B_{R}(p_{\infty})\subset C(X)$ such that $u(p_{\infty})=0$.
We define functions $I_{u}(s)$, $D_{u}(s)$, $U_{u}(s)$ and $J_{u}(s)$ for $s\in(0,R)$ as follows:
\begin{align}
I_{u}(s):=\int_{X}u^{2}(x,s)dm_{X},
\end{align}
\begin{align}
D_{u}(s):=s^{2-n}\int_{B_{s}(p_{\infty})}|Du|^{2}dm_{C(X)},
\end{align}
\begin{align}
U_{u}(s)=\frac{D_{u}(s)}{I_{u}(s)},
\end{align}
\begin{align}
J_{u}(s)=\frac{1}{m_{C(X)}(B_{s}(p_{\infty}))}\int_{B_{s}(p_{\infty})}u^{2}dm_{C(X)}.
\end{align}

The definitions of these functions are motivated by \cite{CM97b}.

It is easy to see that
\begin{align}
J_{u}(s)=\frac{1}{m_{C(X)}(B_{s}(p_{\infty}))}\int_{0}^{s}I_{u}(r)r^{n-1}dr.
\end{align}

Making use of Theorem \ref{eigenfcn}, the following proposition is proved in \cite{Huang19} (a detailed proof of some facts in the following can also be found in Proposition 3.2 in \cite{Huang20}).
\begin{prop}\label{prop3.3}
Suppose $u(x,r)=\sum_{i=1}^{\infty}c_{i}r^{\alpha_{i}}\varphi_{i}(x)$ is a non-constant harmonic function on $B_{R}(p_{\infty})$ such that $u(p_{\infty})=0$ as in Theorem \ref{eigenfcn}, then the following equalities hold for every $s, r\in(0,R)$:
\begin{align}\label{3.1111}
I_{u}(s)=\sum_{i=1}^{\infty}c_{i}^{2}s^{2\alpha_{i}},
\end{align}
\begin{align}\label{3.2222}
D_{u}(s)=\sum_{i=1}^{\infty}c_{i}^{2}\alpha_{i}s^{2\alpha_{i}},
\end{align}
\begin{align}\label{3.5555}
I_{u}(s)=I_{u}(r)\exp\biggl(\int_{r}^{s}\frac{2U_{u}(t)}{t}dt\biggr).
\end{align}
Furthermore, $U_{u}(s)$ is a non-decreasing function for $s$.
\end{prop}

Given $k>0$, denote by
$$\mathcal{H}_{k}(C(X)):=\{u\in W^{1,2}_{\mathrm{loc}}(M)\mid \Delta u=0, |u(y)|\leq C(d_{C(X)}(y,p_{\infty})+1)^{k} \text{ for some }C\}.$$

\begin{prop}\label{prop3.4}
For any $k>0$, we have
\begin{align}\label{3.1222}
\mathrm{dim}\mathcal{H}_{k}(C(X))=N_{X}(k(k+n-2)).
\end{align}
For any $\beta>k$ such that $(k,\beta]\cap\mathscr{D}_{X}=\emptyset$, we have $\mathcal{H}_{\beta}(C(X))=\mathcal{H}_{k}(C(X))$.
In addition, for any $u\in \mathcal{H}_{k}(C(X))$ with $u(p_{\infty})=0$,
\begin{align}\label{3.1234}
J_{u}(s)
\leq 2^{2\alpha_{N}}J_{u}(\frac{s}{2})
\end{align}
holds for every $s>0$.
In (\ref{3.1234}) we denote by $N=N_{X}(k(k+n-2))$ for simplicity, and $\alpha_{N}\geq0$ is given by  $\lambda_{N}=\alpha_{N}(\alpha_{N}+n-2)$, where $\lambda_{N}$ is the $N$-th eigenvalue of $(X,d_{X},m_{X})$.
\end{prop}
\begin{proof}
For any nonzero $u\in \mathcal{H}_{k}(C(X))$ with $u(p_{\infty})=0$, for any $s>0$, we have
\begin{align}\label{3.1112}
I_{u}(s)=\int_{X}u(x,s)^{2}dm_{X}\leq C(1+s^{2k})
\end{align}
for some $C>0$.
By Theorem \ref{eigenfcn}, we write $u(x,r)=\sum_{i=1}^{\infty}c_{i}r^{\alpha_{i}}\varphi_{i}(x)$.
By (\ref{3.1112}) and (\ref{3.1111}),
\begin{align}
\sum_{i=1}^{\infty}c_{i}^{2}s^{2\alpha_{i}}\leq C(1+s^{2k})
\end{align}
holds for any $s>0$.
Hence $c_{i}=0$ for any $i$ such that $\alpha_{i}> k$.
Thus $\mathcal{H}_{k}(C(X))$ is spanned by those $r^{\alpha_{i}}\varphi_{i}(x)$ with $0\leq i\leq N=N_{X}(k(k+n-2))$.
This proves (\ref{3.1222}).

For any $\beta>k$ such that $(k,\beta]\cap\mathscr{D}_{X}=\emptyset$, we have $N_{X}(\beta(\beta+n-2))=N_{X}(k(k+n-2))$, and hence
$\mathcal{H}_{\beta}(C(X))=\mathcal{H}_{k}(C(X))$.

Furthermore, since
\begin{align}
J_{u}(s)=\sum_{i=1}^{N}\frac{c_{i}^{2}}{2\alpha_{i}+n}s^{2\alpha_{i}},
\end{align}
we have
\begin{align}
\frac{J_{u}(s)}{J_{u}(\frac{s}{2})}
=\frac{\sum_{i=1}^{N}\frac{c_{i}^{2}}{2\alpha_{i}+n}s^{2\alpha_{i}}}{\sum_{i=1}^{N} \frac{c_{i}^{2}}{2\alpha_{i}+n}(\frac{s}{2})^{2\alpha_{i}}}
\leq 2^{2\alpha_{N}}
\end{align}
for every $s>0$.
\end{proof}

\section{Proof of main theorems}\label{sec-4}

In this section we give the proof of the main theorems.
Let $(M^{n},g)$ be a Riemannian manifold with nonnegative Ricci curvature, $\mu$ be the volume measure determined by the metric $g$, $p$ be a fixed point on $M$.
Given a function $u$ defined on $B_{R}(p)\subset M$, define $J_{u}:(0,R)\rightarrow\mathbb{R}$ to be
\begin{align}
J_{u}(r)=\bbint_{B_{r}(p)}u^{2}d\mu.
\end{align}

Firstly we recall the following theorem, which is a special case of the three circles theorem proved by Xu:

\begin{thm}[see Theorem 3.2 in \cite{X15}]\label{three-circles}
Suppose $(M^{n},g)$ satisfies the assumptions in Theorem \ref{main-5}, with the unique tangent cone at infinity, denoted by $(C(X),d_{C(X)},m_{C(X)})$.
Let $\mathscr{D}_{X}$ be given by (\ref{Dx}).
Then for any $\beta\notin \mathscr{D}_{X}$, there exists $T=T(\beta, M)>1$ such that for every $r\geq T$ the following statement is true.
Suppose $u(x)$ is any harmonic function defined on $B_{r}(p)$, and if
$$J_{u}(r)\leq2^{2\beta}J_{u}(\frac{r}{2})$$
holds, then
$$J_{u}(\frac{r}{2})\leq2^{2\beta}J_{u}(\frac{r}{4}).$$
\end{thm}

The following is a three circles theorem for harmonic functions in $\mathcal{H}_{k}(M)$.
The proof is mainly modified from that of Theorem 1.4 in \cite{Huang19}.

\begin{prop}\label{main-8-weak}
Suppose $(M^{n},g)$ satisfies the assumptions in Theorem \ref{main-5}, with the unique tangent cone at infinity, denoted by $(C(X),d_{C(X)},m_{C(X)})$.
Let $\mathscr{D}_{X}$ be given by (\ref{Dx}).
Given any $k>0$ with $k(k+n-2)\geq\lambda_{1}$, the following holds.
For any $\epsilon>0$ sufficiently small such that $(k,k+2\epsilon)\cap\mathscr{D}_{X}=\emptyset$, there exists a $T=T(M,k+\epsilon)$ such that for any $u \in\mathcal{H}_{k}(M)$ with $u(p)=0$,
\begin{description}
  \item[(1)] for any $r\geq T$, we have
  \begin{align}\label{5.33333123456}
  J_{u}(2r) \leq 2^{2k+2\epsilon}J_{u}(r);
  \end{align}
  \item[(2)] for any $s\geq t\geq T$, we have
  \begin{align}\label{5.00113}
  J_{u}(s)\leq C\bigl(\frac{s}{t}\bigr)^{2k+2\epsilon} J_{u}(t),
  \end{align}
  where $C$ is a constant depending only on $n$ and $k$.
\end{description}
\end{prop}

\begin{proof} [Proof of Proposition \ref{main-8-weak}]

Given a non-trivial harmonic function $u\in\mathcal{H}_{k}(M)$ such that $u(p)=0$.
By the volume doubling property, it is easy to check that (2) can be derived from (1).
Hence we only need to prove (1).

Firstly, we define a function $f:[0,\infty)\rightarrow[0,\infty)$ by
\begin{align}\label{4.1}
f(s)=\int_{B_{s}(p)}|\nabla u|^{2}d\mu.
\end{align}
Then it is not hard to check that $f$ satisfies the following properties:
\begin{description}
  \item[(a)] $f$ is a non-decreasing function;
  \item[(b)] $f(s)>0$ for $s>0$;
  \item[(c)] $f(s)\leq C(1+s^{2k+n-2})$.
\end{description}
Here and in the following, $C$ will denote constants depending only on $n$ and $u$, and $C$ may be different in different lines.

Applying Lemma 3.1 in \cite{CM97a} to $f$, for any $N\geq2$, there exists a sequence of positive integers $\{m_{i}\}$ (depending on $N$) such that $m_{i}\rightarrow\infty$ and
\begin{align}\label{5.1111}
f(N^{m_{i}+1})\leq 2N^{2k+n-2}f(N^{m_{i}}).
\end{align}

Now suppose $N\gg1$ is a fixed number.

By the mean value inequality (see \cite{LS84}), we have
\begin{align}\label{5.2}
\sup_{B_{\frac{1}{5}N^{m_{i}+1}}(p)}|\nabla u|^{2}\leq C\frac{f(\frac{1}{4}N^{m_{i}+1})}{\mu(B_{\frac{1}{4}N^{m_{i}+1}}(p))}.
\end{align}

Denote by $g_{i}=(N^{m_{i}})^{-2}g$. Suppose $\rho_{i}$ is the distance determined by $g_{i}$, $\mu_{i}$ is the volume measure determined by $g_{i}$, and $\nu_{i}$ is the renormalized measure.
Suppose up to passing to a subsequence, $(M_{i},p_{i},\rho_{i},\nu_{i})$ converge in the pointed measured Gromov-Hausdorff sense to the tangent cone at infinity $(C(X),p_{\infty},d_{C(X)},m_{C(X)})$,
where $M_{i}$ is the same manifold as $M$, $p_{i}=p$.
Denote by $\nabla^{(i)}$ the gradient with respect to $g_{i}$, by $D$ the Cheeger derivative with respect to $(C(X),d_{C(X)},m_{C(X)})$.
Denote by $B_{r}(p_{i})=\{x\in M_{i}\mid\rho_{i}(x,p_{i})< r\}$, $B_{r}(p_{\infty})=\{x\in C(X)\mid d_{C(X)}(x,p_{\infty})< r\}$.

Define
\begin{align}\label{ui}
u_{i}(x)=\frac{u(x)}{N^{m_{i}}\bigl(\frac{f(N^{m_{i}})}{\mu(B_{N^{m_{i}}}(p))}\bigr)^{\frac{1}{2}}},
\end{align}
and denote by
\begin{align}
J_{u_{i}}^{(i)}(s)=\bbint_{B_{s}(p_{i})}u_{i}^{2}d\mu_{i}.
\end{align}
Obviously, $u_{i}$ is a harmonic function on $M_{i}$ and $u_{i}(p_{i})=0$.
By (\ref{5.1111}) and (\ref{5.2}), we have
\begin{align}
\sup_{B_{\frac{1}{5}N}(p_{i})}|\nabla^{(i)} u_{i}|^{2}
=& \frac{\mu(B_{N^{m_{i}}}(p))}{f(N^{m_{i}})} \sup_{B_{\frac{1}{5}N^{m_{i}+1}}(p)}|\nabla u|^{2} \\
\leq& C\frac{\mu(B_{N^{m_{i}}}(p))}{\mu(B_{\frac{1}{4}N^{m_{i}+1}}(p))} \cdot\frac{f(\frac{1}{4}N^{m_{i}+1})}{f(N^{m_{i}})}\nonumber\\
\leq& CN^{2k+n-2}.\nonumber
\end{align}
Hence each $u_{i}$ is a Lipschitz harmonic function on $B_{\frac{1}{10}N}(p_{i})$ with uniform Lipschitz constant $CN^{k+\frac{n}{2}-1}$, and by Propositions \ref{AA} and \ref{2.7777}, up to passing to a subsequence, $u_{i}$ converge uniformly to some Lipschitz harmonic function $\omega_{N}$ defined on $B_{\frac{1}{10}N}(p_{\infty})\subset C(X)$ with $\omega_{N}(p_{\infty})=0$.

We define functions $I_{\omega_{N}}(s)$, $D_{\omega_{N}}(s)$, $U_{\omega_{N}}(s)$ and $J_{\omega_{N}}(s)$, where $s\in(0,\frac{1}{10}N)$, as in Section \ref{sec-3}.

Note that for any $s\in(0,\frac{1}{10}N)$, it is easy to see
\begin{align}\label{5.21111}
\int_{B_{s}(p_{i})}|\nabla^{(i)}u_{i}|^{2}d\nu_{i}
=&\frac{f(sN^{m_{i}})}{f(N^{m_{i}})}.
\end{align}

By Proposition \ref{2.7777}, we have
$$\int_{B_{1}(p_{\infty})}|D\omega_{N}|^{2}dm_{C(X)}=\lim_{i\rightarrow\infty}\int_{B_{1}(p_{i})} |\nabla^{(i)}u_{i}|^{2}d\nu_{i}=1,$$
i.e. $D_{\omega_{N}}(1)=1$.

Furthermore,
\begin{align}\label{5.5555}
\frac{I_{\omega_{N}}(\frac{N}{100})}{I_{\omega_{N}}(1)}
=\frac{D_{\omega_{N}}(\frac{N}{100})}{D_{\omega_{N}}(1)}\frac{U_{\omega_{N}}(1)}{U_{\omega_{N}}(\frac{N}{100})}
\leq D_{\omega_{N}}(\frac{N}{100}),
\end{align}
where in the last equality we use the fact that $U_{\omega_{N}}(s)$ is non-decreasing, and we assume $N>100$.

On the other hand, by Proposition \ref{2.7777}, (\ref{5.21111}) and (\ref{5.1111}), we have
\begin{align}\label{5.6666}
D_{\omega_{N}}(\frac{N}{100})
&=\biggl(\frac{N}{100}\biggr)^{2-n}\int_{B_{\frac{N}{100}}(p_{\infty})}|D\omega_{N}|^{2}dm_{C(X)}\\
&{=}\biggl(\frac{N}{100}\biggr)^{2-n}\lim_{i\rightarrow\infty}
\int_{B_{\frac{N}{100}}(p_{i})}|\nabla^{(i)}u_{i}|^{2}d\nu_{i}\nonumber\\
&{\leq}\biggl(\frac{N}{100}\biggr)^{2-n}
\limsup_{i\rightarrow\infty}\frac{f(N^{m_{i}+1})}{f(N^{m_{i}})}
\nonumber\\
&\leq CN^{2k}.\nonumber
\end{align}

Combining (\ref{5.5555}), (\ref{5.6666}) and (\ref{3.5555}), we have
\begin{align}
\exp\biggl(\int_{1}^{\frac{N}{100}}\frac{2U_{\omega_{N}}(t)}{t}dt\biggr)\leq CN^{2k}.
\end{align}

Using the non-decreasing property of $U_{\omega_{N}}(s)$ again, for any $L\in[1,\frac{N}{100})$, we have
\begin{align}
CN^{2k}\geq\exp\biggl(\int_{L}^{\frac{N}{100}}\frac{2U_{\omega_{N}}(L)}{t}\biggr)
=\biggl(\frac{N}{100L}\biggr)^{2U_{\omega_{N}}(L)},
\end{align}
and then
\begin{align}\label{4.1234}
2U_{\omega_{N}}(L)\leq\frac{\log C+2k\log N}{\log N-\log(100L)}.
\end{align}

Thus, given any $\epsilon>0$, for all $L\geq 1$, there exists a constant $\tilde{C}=\tilde{C}(\epsilon,L,n,u)$ such that for any $N\geq \tilde{C}$, we have
\begin{align}
U_{\omega_{N}}(s)\leq U_{\omega_{N}}(L)<k+\frac{\epsilon}{2}, \qquad \forall s\in(0,L].
\end{align}
Thus we have
\begin{align}\label{5.9999}
I_{\omega_{N}}(s)=I_{\omega_{N}}(\frac{s}{2})\exp\biggl(\int_{\frac{s}{2}}^{s}\frac{2U_{\omega_{N}}(t)}{t}dt \biggr)\leq 2^{2k+\epsilon}I_{\omega_{N}}(\frac{s}{2}), \qquad \forall s\in(0,L].
\end{align}
Then for every $s\in(0,L]$, we have
\begin{align}\label{5.999999}
J_{\omega_{N}}(s)&=\frac{1}{m_{C(X)}(B_{s}(p))}\int_{0}^{s}I_{\omega_{N}}(r)r^{n-1}dr\\
&\leq\frac{1}{m_{C(X)}(B_{s}(p))}\int_{0}^{s}I_{\omega_{N}}(\frac{r}{2})2^{2k+\epsilon}r^{n-1}dr \nonumber\\
&= 2^{2k+\epsilon} \frac{1}{m_{C(X)}(B_{\frac{s}{2}}(p))}\int_{0}^{\frac{s}{2}}I_{\omega_{N}}(r)r^{n-1}dr \nonumber\\
&= 2^{2k+\epsilon}J_{\omega_{N}}(\frac{s}{2}).\nonumber
\end{align}

We will take $\epsilon$ to be any sufficiently small positive number such that $(k,k+2\epsilon)\cap\mathscr{D}(M)=\emptyset$.
In the following, we denote by $\beta=k+\epsilon$ for simplicity.
Once $\epsilon$ is fixed, we can fix $N$ such that
(\ref{5.999999}) holds for $s\in(0,2]$.

Denote by $r_{i}=N^{m_{i}}$.
By Proposition \ref{2.7777} and (\ref{5.999999}), for any $s\in[\frac{1}{2},1]$, we have
\begin{align}\label{5.0011}
\lim_{i\rightarrow\infty}\frac{J_{u}(sr_{i})}{J_{u}(\frac{s}{2}r_{i})} =\lim_{i\rightarrow\infty}\frac{J^{(i)}_{u_{i}}(s)}{J^{(i)}_{u_{i}}(\frac{s}{2})}
\overset{(\ast)}{=}\frac{J_{\omega_{N}}(s)}{J_{\omega_{N}}(\frac{s}{2})}\leq2^{2k+\epsilon}.
\end{align}
It is not hard to check that the convergence in $(\ast)$ is uniformly on $s\in [\frac{1}{2},1]$.
Hence for sufficiently large $i$, and any $s\in [\frac{1}{2},1]$, we have
\begin{align}\label{5.00111}
J_{u}(sr_{i})\leq 2^{2\beta} J_{u}(\frac{s}{2}r_{i}).
\end{align}

By Theorem \ref{three-circles}, we derive that
\begin{align}\label{5.00112}
J_{u}(s)\leq 2^{2\beta} J_{u}(\frac{s}{2})
\end{align}
holds for every $r_{i}\geq s\geq T$, where $T=T(\beta, M)$ is given in Theorem \ref{three-circles}.
Since $r_{i}\rightarrow\infty$, we obtain that (\ref{5.00112}) holds for every $s\geq T$.
The proof is completed.
\end{proof}

Now we can prove Theorem \ref{main-5}.

\begin{proof}[Proof of Theorem \ref{main-5}]
Given any $k>0$, we choose $\epsilon>0$ sufficiently small such that $(k,k+2\epsilon)\cap\mathscr{D}_{X}=\emptyset$.
Denote by $W=\{u\in\mathcal{H}_{k}(M)|u(p)=0\}$, $l=h_{k}-1$ and $\beta=k+\epsilon$ for simplicity.
Note that
\begin{align}\label{4.18}
N_{X}(k(k+n-2))=N_{X}(\beta(\beta+n-2)).
\end{align}

By Proposition \ref{main-8-weak}, there exists $T=T(M,\beta)>0$ and $C=C(n,\beta)$ such that
\begin{align}\label{5.1234}
J_{u}(s)\leq C\bigl(\frac{s}{t}\bigr)^{2\beta} J_{u}(t)
\end{align}
holds for any $u\in W$ and any $s\geq t\geq T$.

Given $r_{i}\rightarrow \infty$ and we consider $g_{i}=r_{i}^{-2}g$.
Suppose up to passing to a subsequence, $(M_{i},p_{i},\rho_{i},\nu_{i})$ converge in the pointed measured Gromov-Hausdorff sense to the tangent cone at infinity $(C(X),p_{\infty},d_{C(X)},m_{C(X)})$,
where $M_{i}$ is the same manifold as $M$, $p_{i}=p$, $\rho_{i}$ is the distance determined by $g_{i}$, $\mu_{i}$ is the volume measure determined by $g_{i}$, and $\nu_{i}$ is the renormalized measure.
Denote by $\nabla^{(i)}$ the gradient with respect to $g_{i}$, by $D$ the Cheeger derivative with respect to $(C(X),d_{C(X)},m_{C(X)})$.
Denote by $B_{r}(p_{i})=\{x\in M_{i}\mid\rho_{i}(x,p_{i})< r\}$, $B_{r}(p_{\infty})=\{x\in C(X)\mid d_{C(X)}(x,p_{\infty})< r\}$,
and denote by $J_{u}^{(i)}(s)=\bbint_{B_{s}(p_{i})}u^{2}d\mu_{i}$ for any function $u: M_{i}\rightarrow\mathbb{R}$.

For $i$ sufficiently large, we can find $\{v_{i,a}\}_{1\leq a\leq l}\subset W$ such that
\begin{align}\label{6.111}
\int_{B_{1}(p_{i})}v_{i,a}v_{i,b} d\nu_{i}=\delta_{ab}
\end{align}
for $a, b\in\{1,\ldots,l\}$.
Then by (\ref{5.1234}),
\begin{align}
J^{(i)}_{v_{i,a}}(s)\leq C(n,\beta) s^{2\beta}
\end{align}
holds for every $a\in\{1,\ldots,l\}$ and $s\geq 1$.

By mean value inequality, for any $a\in\{1,\ldots,l\}$ and $s\geq1$, we have
\begin{align}
\sup_{B_{s}(p_{i})} |v_{i,a}|^{2}\leq C(n)J^{(i)}_{v_{i,a}}(2s)\leq C(n,\beta) s^{2\beta},
\end{align}
and then by Cheng-Yau's gradient estimate (\cite{CY75}), we have
\begin{align}
\sup_{B_{s}(p_{i})} |\nabla^{(i)}v_{i,a}|^{2}\leq  C(n,\beta) s^{2\beta-2}.
\end{align}
Thus by Propositions \ref{AA} and \ref{2.7777}, for each $a\in\{1,\ldots,l\}$, $v_{i,a}$ converge locally uniformly to a harmonic function $v_{a}:C(X)\rightarrow \mathbb{R}$ with $v_{a}(p_{\infty})=0$, and
\begin{align}
\sup_{B_{s}(p_{\infty})} |v_{a}|\leq  C(n,\beta) s^{\beta}
\end{align}
holds for any $s\geq1$.
Thus each $v_{a}\in \mathcal{H}_{\beta}(C(X))$.

In addition, by Proposition \ref{2.7777} and (\ref{6.111}), we have
\begin{align}
\int_{B_{1}(p_{\infty})}v_{a}v_{b} d m_{C(X)}=\delta_{ab}
\end{align}
for $a, b\in\{1,\ldots,l\}$.
Thus $v_{1},\ldots,v_{l}$ are linearly independent.
Thus
\begin{align}\label{4.19}
\textmd{dim}\mathcal{H}_{\beta}(C(X))\geq l+1.
\end{align}
Then (\ref{5.33333123}) follows from Proposition \ref{prop3.4}, (\ref{4.18}) and (\ref{4.19}).

The proof is completed.
\end{proof}

Now we prove some corollaries of Theorem \ref{main-5}.

\begin{proof}[Proof of Corollary \ref{main-6}]
Given any $k>0$ such that $k\notin \mathscr{D}_{X}$, by the discreteness of $\mathscr{D}_{X}$, we can find $\beta<k$ such that
\begin{align}\label{4.27}
N_{X}(\beta(\beta+n-2))=N_{X}(k(k+n-2)).
\end{align}
Then (\ref{5.33333123ggggg}) follows from (\ref{4.27}), (\ref{5.33333}) and (\ref{5.33333123}).
\end{proof}

\begin{proof}[Proof of Corollary \ref{main-7}]
By Corollary \ref{main-6} and Proposition \ref{prop2.1111}, we have
\begin{align}
\lim_{k\rightarrow\infty}\frac{h_{k}}{k^{n-1}}=\lim_{k\rightarrow\infty}\frac{N_{X}(k(k+n-2))}{k^{n-1}}=\frac{n\omega_{n-1}\alpha}{(2\pi)^{n-1}}
=\frac{2\alpha}{(n-1)!\omega_{n}}.
\end{align}
\end{proof}

Finally we prove Theorem \ref{main-8} in the following.

\begin{proof}[Proof of Theorem \ref{main-8}]

Given any positive sequence $\{r_{i}\}$ with $r_{i}\rightarrow \infty$, we consider $g_{i}=r_{i}^{-2}g$.
Suppose up to passing to a subsequence, $(M_{i},p_{i},\rho_{i},\nu_{i})$ converge in the pointed measured Gromov-Hausdorff sense to the tangent cone at infinity $(C(X),p_{\infty},d_{C(X)},m_{C(X)})$, where we use the same notations as in the proof of Theorem \ref{main-5}.

Suppose $u_{i}:M_{i}\rightarrow\mathbb{R}$ is given by $u_{i}(x)=\frac{1}{(J_{u}(r_{i}))^{\frac{1}{2}}}u(x)$.
Then
\begin{align}
\int_{B_{1}(p_{i})}u_{i}^{2} d\nu_{i}=1.
\end{align}

Given any $\delta>0$ sufficiently small such that $(k,k+2\delta)\cap\mathscr{D}_{X}=\emptyset$, denote by $\beta=k+\delta$.
By Proposition \ref{main-8-weak}, we derive that for any sufficiently large $i$,
\begin{align}
J_{u_{i}}^{(i)}(s)\leq C(n,k)s^{2\beta}
\end{align}
holds for any $s\geq 1$.
Hence by the same argument as in the proof of Theorem \ref{main-5}, we know $u_{i}$ converges locally uniformly to a harmonic function $v\in \mathcal{H}_{\beta}(C(X))$ with $v(p_{\infty})=0$.

By Proposition \ref{prop3.4}, we have $v(x,r)=\sum_{i=1}^{N}c_{i}r^{\alpha_{i}}\varphi_{i}(x)$, where $N=N_{X}(\beta(\beta+n-2))=N_{X}(k(k+n-2))$.
Let $P(v)$ denote the smallest integer $a$ in $\{1,\ldots,N\}$ such that $c_{i}=0$ for any $i>a$.
Denote by $Q(v)=\alpha_{P(v)}\in\mathscr{D}_{X}$.

Let $\gamma$ be the smallest one in $\mathscr{D}_{X}$ of the form $Q(v)$, where $v$ is obtained by taking limits by choice of parameters $r_{i}$ and the subsequences as above.

For simplicity of notation, we assume the above $v$ satisfies $Q(v)=\gamma$.
In particular, $v\in \mathcal{H}_{\gamma}(C(X))$.
Hence by Proposition \ref{prop3.4}, for any $s>0$ we have
\begin{align}
J_{v}(s)\leq 2^{2\gamma}J_{v}(\frac{s}{2}).
\end{align}

Following the argument in the last part of the proof of Proposition \ref{main-8-weak}, for any $\epsilon>0$ such that $(\gamma,\gamma+2\epsilon)\cap\mathscr{D}_{X}=\emptyset$,
\begin{align}
J_{u}(s)\leq 2^{2\gamma+2\epsilon}J_{u}(\frac{s}{2})
\end{align}
holds for any $s\geq T=T(M,u,\epsilon)$.
This completes the proof of (1).

Once we have (1), then by mean value inequality and Cheng-Yau's gradient estimate, we derive $u\in\mathcal{H}_{\gamma+\epsilon}(M)$.

Finally, suppose $u\in\mathcal{H}_{\gamma-\epsilon}(M)$ holds for some $\epsilon>0$, then after blowing down the manifold by suitable parameters, we obtain a limit harmonic function $\tilde{v}:C(X)\rightarrow \mathbb{R}$ with $Q(\tilde{v})\leq \gamma-\frac{\epsilon}{2}$, which contradicts to the definition of $\gamma$.
Thus the proof of (3) is completed.
\end{proof}

\section{Remarks on the collapsed cases}\label{sec-5}

In this section, we give some remarks on the estimates of $h_{k}(M)$ when $M$ is collapsed.

As mentioned in Remark \ref{rem1.14}, following the arguments in this paper, we can derive that conclusions similar to Theorem \ref{main-5}, Corollary \ref{main-6}, and Theorem \ref{main-8} still hold on a manifold $(M^{n},g)$ with nonnegative Ricci curvature and satisfying assumptions (1)-(3) in Remark \ref{rem1.14}.
We omit the detailed proofs.
Recall that in \cite{CM98b}, Colding and Minicozzi asked, under some assumptions on the tangent cones at infinity on a collapsed manifold $M$, are there some improved upper of $h_{k}(M)$ similar but better than (\ref{cm_asy})?
See Question 6.53 of \cite{CM98b} for details.
Even though assumptions (1)-(3) in Remark \ref{rem1.14} are more restrictive than the assumptions in Question 6.53 of \cite{CM98b}, our conclusions seems better.

In Theorem 6.58 of \cite{CM98b}, Colding and Minicozzi proved an improved upper of $h_{k}(M)$ (see (6.59) in \cite{CM98b}) provided $(M^{n},g)$ have nonnegative sectional curvature and its tangent cone at infinity has Hausdorff dimension $m$.
Recall that if $(M^{n},g)$ have nonnegative sectional curvature, then no matter $M$ has maximal volume growth, its tangent cone at infinity, denote by $M_{\infty}$, is unique.
Furthermore, $M_{\infty}$ is an Alexandrov space with nonnegative curvature with integral Hausdorff dimension $m\leq n$; and if $m\geq2$, then $M_{\infty}$ is an Euclidean cone over an $(m-1)$-dimensional Alexandrov space $X$ with curvature bounded from below by $1$; if $m=1$, then either $M_{\infty}$ is isometric to a line or a half line.
Note that if $m<n$, then in general the renormlized limit measure on $M_{\infty}$ depends on the parameters to blow down $M$, see Example 1.24 in \cite{CC97}.

If $M_{\infty}$ is isometric to a line, then $M$ itself splits.
In the case that $M_{\infty}$ is isometric to a half line, then argue as in the proof of Theorem 1.4 of \cite{Huang20II}, we derive that there is no non-constant harmonic function with polynomial growth on $M$.

Thus in the following we always assume $m\geq2$, in particular assumption (2) in Remark \ref{rem1.14} holds.

In the following, we always assume:
\begin{description}
  \item[(1)'] the renormlized limit measures on $M_{\infty}$ is unique, and up to a constant, it is equal to the $m$-dimensional Hausdorff measure of $M_{\infty}$.
\end{description}
Hence $M$ satisfies assumption (1) of Remark \ref{rem1.14} with $\kappa=m$.
In this case, by Theorem 1.2 in \cite{Huang19},
\begin{align}\label{5.333338652}
h_{k}\geq N_{X}(\beta(\beta+m-2))
\end{align}
holds for any positive number $\beta<k$.

On the other hand, if we denote by $V=\mathcal{H}^{m-1}(X)$ (hence the $m$-dimensional Hausdorff measure of the unit ball centered at the tips in $M_{\infty}$ equal to $\frac{1}{m}V$), then applying Weyl's law (see Corollary 4.8 of \cite{AHT17}) to the Alexandrov space $X$, we have
\begin{align}\label{1.1111789}
\lim_{\lambda\rightarrow\infty}\frac{N_{X}(\lambda)}{\lambda^{\frac{m-1}{2}}}
=\frac{\omega_{m-1}V}{(2\pi)^{m-1}}=\frac{2V}{m!\omega_{m}}.
\end{align}

By (\ref{5.333338652}) and (\ref{1.1111789}), we can derive
\begin{align}\label{5.6655}
\liminf_{k\rightarrow\infty}k^{1-m}h_{k}
\geq \frac{2V}{m!\omega_{m}}.
\end{align}

We want to estimate the upper bound of $h_{k}$, following the argument in the proof of Theorem \ref{main-5}.
Note that we don't require $M$ to satisfy assumption (3) in Remark \ref{rem1.14}.
By volume comparision, we have the following natural volume upper bound:
\begin{align}
\frac{\mu(B_{r}(p))}{r^{n}}\leq \omega_{n}, \qquad \forall r>0.
\end{align}
We follow the argument in the proof of Proposition \ref{main-8-weak}, define $f$ as in (\ref{4.1}) for $u\in\mathcal{H}_{k}(M)$ with $u(p)=0$.
Then (\ref{5.1111}) still holds.
We consider $g_{i}=(N^{m_{i}})^{-2}g$, and use the same notations as in the proof of Proposition \ref{main-8-weak}, and the function $u_{i}$ defined in (\ref{ui}) converges to a harmonic function $\omega_{N}$ defined on $B_{\frac{1}{10}N}(p_{\infty})\subset C(X)=M_{\infty}$.
By the conic structure, $\omega_{N}$ has a beautiful form, see Theorem 3.1 in \cite{Huang19}.
Then we define functions $I_{\omega_{N}}(s)$, $D_{\omega_{N}}(s)$, $U_{\omega_{N}}(s)$ and $J_{\omega_{N}}(s)$.
We remark that in this case, $D_{\omega_{N}}:(0,\frac{1}{10}N)\rightarrow \mathbb{R}$ is defined to be
\begin{align}
D_{\omega_{N}}(s):=s^{2-m}\int_{B_{s}(p_{\infty})}|D\omega_{N}|^{2}dm_{C(X)}.
\end{align}
Because of this difference, we have
\begin{align}
D_{\omega_{N}}(\frac{N}{100})&{\leq}\biggl(\frac{N}{100}\biggr)^{2-m}
\limsup_{i\rightarrow\infty}\frac{f(N^{m_{i}+1})}{f(N^{m_{i}})}
\leq CN^{2k+n-m}
\end{align}
instead of (\ref{5.6666}).
Then we derive
\begin{align}
2U_{\omega_{N}}(L)\leq\frac{\log C+(2k+n-m)\log N}{\log N-\log(100L)}.
\end{align}
instead of (\ref{4.1234}).
Thus, given any $\epsilon>0$, for all $L\geq 1$, there exists a constant $\tilde{C}=\tilde{C}(\epsilon,L,n,u)$ such that for any $N\geq \tilde{C}$, we have
\begin{align}
U_{\omega_{N}}(s)\leq U_{\omega_{N}}(L)<k+\frac{n-m}{2}+\frac{\epsilon}{2}, \qquad \forall s\in(0,L].
\end{align}
Then for every $s\in(0,L]$, we have
\begin{align}\label{5.123333}
J_{\omega_{N}}(s)\leq 2^{2k+n-m+\epsilon}J_{\omega_{N}}(\frac{s}{2})
\end{align}
instead of (\ref{5.999999}).
We take $\epsilon$ to be any sufficiently small positive number such that $(k+\frac{n-m}{2},k+\frac{n-m}{2}+2\epsilon)\cap\mathscr{D}(M)=\emptyset$.
Combing Proposition \ref{2.7777}, (\ref{5.123333}) and Xu's three circle theorem (Theorem 3.2 in \cite{X15}),
we derive that
\begin{align}\label{5.0011277}
J_{u}(s)\leq 2^{2k+n-m+2\epsilon} J_{u}(\frac{s}{2})
\end{align}
holds for every $s\geq T$, where $T=T(k+\frac{n-m}{2}+\epsilon, M)$ is given in Theorem 3.2 in \cite{X15}.
With the help of (\ref{5.0011277}), we follow the argument in the proof of Theorem \ref{main-5} to conclude that
\begin{align}\label{5.222}
h_{k}\leq N_{X}((k+\frac{n-m}{2})(k+\frac{n+m}{2}-2)).
\end{align}
Combing (\ref{1.1111789}) and (\ref{5.222}) we have
\begin{align}
\limsup_{k\rightarrow\infty}k^{1-m}h_{k}
\leq&\lim_{k\rightarrow\infty}\frac{N_{X}((k+\frac{n-m}{2})(k+\frac{n+m}{2}-2))}{k^{m-1}} \\ =&\frac{2V}{m!\omega_{m}}.\nonumber
\end{align}

In conclusion, if $(M^{n},g)$ has nonnegative sectional curvature and satisfies the assumption (1)' with $m\geq2$, then we have
\begin{align}
\lim_{k\rightarrow\infty}k^{1-m}h_{k}
=\frac{2V}{m!\omega_{m}}.
\end{align}
This conclusion is better than (6.59) of \cite{CM98b}.

\end{document}